\documentclass[10pt,article]{amsart}
\usepackage{amssymb}
\begin{document}
\begin{center}
{\bf Euler type integral operator involving k-Mittag-Leffler function }
\end{center}
\begin{center}
\vspace{.2cm} {\bf Waseem A. Khan$^{1}$, Nisar K S$^{2,*}$ and Moin Ahmad$^{3}$}
\end{center}
\begin{center}
{$^{1,3}$Department of Mathematics, Faculty of Science, Integral University,
Lucknow-226026, India}
\end{center}
\begin{center}
{$^{2,*}$Department of Mathematics, College of Arts and Science at Wadi Aldawaser, 11991,
Prince Sattam bin Abdulaziz University, Alkharj, Kingdom of Saudi Arabia}
\end{center}

\begin{center}
{E-mail: waseem08\_khan@rediffmail.com, n.sooppy@psau.edu.sa, moinah1986@gmail.com}
\end{center}

\vspace{.35cm}
\parindent=0mm

\noindent{\bf Abstract.} This paper deals with a Euler type integral operator involving k-Mittag-Leffler function defined by Gupta and Parihar [8]. Furthermore, some special cases are also taken into consideration.

\noindent{\bf Keywords:} Euler type integrals, extended k-beta function, generalized k-Mittag-Leffler function, generalized k-Wright function.

{\bf 2010 Mathematics Subject Classification.:} 33E12, 33C45.

\vspace{.35cm}
\parindent=0mm
\noindent
\par

\noindent{\bf 1. Introduction }

\noindent\par
Many authors namely, Diaz et al. [6], Kokologiannaki [11], Krasniqi [12], Mansour [17], Merovci [15], had introduced k-generalized gamma, k-Zeta, k-Beta functions. They had proven a number of their properties and inequalities for the above k-generalization functions. They had studied k-hypergeometric functions based on k-Pochhamer symbols  for factorial functions. Propose our present research, we begin by mentioning the following definitions of some well known functions:

The integral representation of the k-gamma function as:
$${\Gamma}_{k}(z)=k^{\frac{z}{k}-1}\Gamma{(\frac{z}{k})}=\int_{0}^{\infty}e^{\frac{-t^{k}}{k}}t^{z-1}dt, k \in\mathbb R,z \in \mathbb C.\eqno(1.1)$$
and k-beta function is defined as:
$$B_{k}(x,y)= \frac{1}{k}\int\limits_{0}^{1}t^{\frac{x}{k}-1}(1-t)^{\frac{y}{k}-1}dt=\frac{\Gamma_{k}(x)\Gamma_{k}(y)}{\Gamma_{k}(x+y)}, x>0,y>0.\eqno(1.2) $$

The generalized k-Wright function [7] represented as follows:
$${}_p\Psi_{q}^{k}\left[\begin{array}{lll}
(\alpha_{1},A_{1}),...,(\alpha_{p},A_{p});\\
&z\\
(\beta_{1},B_{1}),...,(\beta_{p},B_{p});\end{array}\right]={}_p\Psi_{q}^{k}\left((\alpha_{j},A_{j})_{1,p};(\beta_{j},B_{j})_{1,q};z\right)$$
$$=\sum\limits_{n=0}^{\infty}\frac{\Gamma_{k}(\alpha_{1}+n A_{1})...,\Gamma_{k}(\alpha_{p}+nA_{p})}{\Gamma_{k}(\beta_{1}+n B_{1})...,\Gamma_{k}(\beta_{p}+n B_{p})}\frac{z^{n}}{n!}.\eqno(1.3)$$

In 1903, the Swedish Mathematician introduced the Mittag-Leffler function $E_{\alpha}(z)$ [16] as:
$$ E_{\alpha}(z) = \sum\limits_{n=0}^{\infty}\frac{z^{n}}{\Gamma(\alpha n +1)}, \eqno(1.4)$$
where $z \in \mathbb C$ and $\Gamma(s)$ is the Gamma function; ${\alpha}\geq 0$.

The Mittag-Leffler function is a direct generalization of $\exp(z)$ in which $\alpha = 1$. Mittag-Leffler function naturally occurs as the solution of fractional order differential equation or fractional order integral equations.

A generalization of $E_{\alpha}(z)$ was studied by Wiman [26, 27] where he defined the function $E_{\alpha,\beta}(z)$ as:

$$E_{\alpha,\beta}(z)= \sum\limits_{n=0}^{\infty}\dfrac{z^{n}}{\Gamma(\alpha n +\beta)},\eqno(1.5)$$

where $\alpha,\beta\in\mathbb C; \Re(\alpha)>0,\Re(\beta)>0.$ Which is also known as Mittag-Leffler function or Wiman's function.

Prabhakar [18] introduced the function $E_{\alpha,\beta}^{\gamma}(z)$ in the form (see  Killbas et al. [9]):

$$E_{\alpha,\beta}^{\gamma}(z)= \sum\limits_{n=0}^{\infty}\dfrac{{(\gamma)_{n}}}{\Gamma(\alpha n +\beta)}\dfrac{z^{n}}{n!},\eqno(1.6)$$

where $\alpha, \beta, \gamma\in\mathbb C  ; \Re(\alpha)>0, \Re(\beta)>0, \Re(\gamma)>0.$

Shukla and Prajapati [24] (see Srivastava and Tomovaski [25]) defined and investigated the function $E_{\alpha,\beta}^{\gamma,q}(z)$ as:

$$E_{\alpha,\beta}^{\gamma,q}(z)= \sum\limits_{n=0}^{\infty}\dfrac{{(\gamma)_{qn}}}{\Gamma(\alpha n +\beta)}\dfrac{z^{n}}{n!},\eqno(1.7)$$

where $\alpha, \beta, \gamma\in\mathbb C  ; \Re(\alpha)>0, \Re(\beta)>0, \Re(\gamma)>0, q \in (0,1)\bigcup \mathbb N $
and $(\gamma)_{qn} =\dfrac{\Gamma(\gamma +qn)}{\Gamma(\gamma)}$ denotes the generalized Pochhammer symbol.

Salim [22] introduced a new generalized Mittag-Leffler function and defined it as:
$$E_{\alpha,\beta}^{\gamma,\delta}(z)= \sum\limits_{n=0}^{\infty}\dfrac{{(\gamma)_{n}}}{\Gamma(\alpha n +\beta)}\dfrac{z^{n}}{(\delta)_{n}}\eqno(1.8)$$

where $\alpha, \beta, \gamma, \delta\in\mathbb C; \Re(\alpha)>0, \Re(\beta)>0, \Re(\gamma)> 0, \Re(\delta)> 0.$

Afterward,  Salim and Faraj [21] introduced the generalized Mittag-Leffler function$E_{\alpha,\beta, p}^{\gamma,\delta,q}(z)$ which is defined as:

$$E_{\alpha,\beta, p}^{\gamma,\delta,q}(z)= \sum\limits_{n=0}^{\infty}\dfrac{{(\gamma)_{qn}}}{\Gamma(\alpha n +\beta)}\dfrac{z^{n}}{(\delta)_{pn}},\eqno(1.9)$$

where ${\alpha,\beta,\gamma,\delta} \in \mathbb C, \min\{\Re(\alpha), \Re(\beta), \Re(\gamma), \Re(\delta)\}>0; p,q>0$ and  ${q < \Re\alpha+p}.$

Currently, Dorrego and Cerutti [5] introduced the generalized k-Mittag-Leffler function as follows:
$$E_{k,\alpha,\beta, p}^{\gamma,\delta,q}(z)= \sum\limits_{n=0}^{\infty}\dfrac{{(\gamma)_{qn,k}}}{{\Gamma}_{k}(\alpha n +\beta)}\dfrac{z^{n}}{(\delta)_{pn,k}},\eqno(1.10)$$

where ${k\in \mathbb R; \alpha,\beta,\gamma,\delta} \in \mathbb C, \min\{\Re(\alpha), \Re(\beta), \Re(\gamma), \Re(\delta)\}>0; p,q>0$ and ${q < \Re\alpha+p}.$

Now, we state the classical beta function denoted by $B(a,b)$ which is defined (see [14], see also [19]) by Euler's integral as:
$$B(a,b)=\int\limits_{0}^{1}t^{a-1}(1-t)^{b-1}dt=\frac{{\Gamma(a)}{\Gamma{(b)}}}{{\Gamma(a+b)}}, (\Re{(a)}>0,\Re(b)>0).\eqno(1.11)$$

In 1997, Chaudhary et al. [4] presented the following extension of Euler's Beta function as follows:
$$B_{\sigma}(a,b)=\int\limits_{0}^{1}t^{a-1}(1-t)^{b-1}\exp\left[-\frac{\sigma}{t(1-t)}\right]dt.\eqno(1.12)$$

In continuation of his work, Lee et al. [13] introduced the generalizations of Euler beta functions and defined it as:
$$B(x,y;p;m)=\int\limits_{0}^{1}t^{x-1}(1-t)^{y-1}\exp\left[-\frac{p}{t^{m}(1-t)^{m}}\right]dt. \eqno(1.13)$$
where $\Re(p) > \Re(m) > 0$.

In this paper, we consider the new generalizations of Euler type k-Beta functions as follows:
$$B_{k}(x,y;p;m)=\int\limits_{0}^{1}t^{\frac{x}{k}-1}(1-t)^{\frac{y}{k}-1}\exp\left[-\frac{p}{t^{\frac{m}{k}}(1-t)^{\frac{m}{k}}}\right]dt,\eqno(1.14)$$
where $k\in\mathbb R;\Re(p) > \Re(m) >0$.

Clearly, when $m= k= 1$, equation (1.14) reduces to (1.12) and further, by taking  $p=0$ in (1.14), we get (1.11).

In this paper, we have obtained some theorems on Euler type integral operator involving generalized k-Mittag-Leffler function and have discussed some special cases.

\noindent{\bf 2. Basic properties of Euler type integral operator involving generalized k-Mittag-Leffler function }
\noindent\par

{\noindent}{\bf Theorem 2.1.} If ${{ k\in \mathbb R}, \alpha,\beta,\gamma,\delta} \in \mathbb C, \Re(\alpha)>0, \Re(\beta)>0, \Re(\gamma)> 0, \Re(\delta)> 0, \Re(A)>0 $, $p,q > 0$ and $q < R{\alpha} + p $, then,
$$\int\limits_{0}^{1}t^{\frac{a}{k}-1}(1-t)^{\frac{b}{k}-1}\exp\left(\frac{-A}{t^{\frac{m}{k}}(1-t)^{\frac{m}{k}}}\right)E_{k,\alpha,\beta, p}^{\gamma,\delta,q}(zt^{\frac{\alpha}{k}})dt$$
$$=\sum\limits_{n=0}^{\infty}\frac{(\gamma)_{qn,k}{z^{n}}}{{\Gamma}_{k}({\alpha}n+\beta)(\delta)_{pn,k}}B_{k}(a+n{\alpha},b;m;A).\eqno(2.1)$$

\noindent{\bf Proof.} In order to derive (2.1), we denote L.H.S. of  (2.1) by $I_1$ and then expanding $E_{k,\alpha,\beta, p}^{\gamma,\delta,q}(zt^{\alpha})$ by using (1.10), to get:

$$I_1=\int\limits_{0}^{1}t^{\frac{a}{k}-1}(1-t)^{\frac{b}{k}-1}\exp\left(\frac{-A}{t^{\frac{m}{k}}(1-t)^{\frac{m}{k}}}\right)\sum\limits_{n=0}^{\infty}\frac{(\gamma)_{qn,k}{z^{n}t^{\frac{n \alpha}{k}}}}{{\Gamma}_{k}({\alpha}n+\beta)(\delta)_{pn,k}}dt$$
Now changing the order of summation and integration (which is guaranteed under the given conditions), we get:
$$I_1=\sum\limits_{n=0}^{\infty}\frac{(\gamma)_{qn,k}{z^{n}}}{{\Gamma}_{k}({\alpha}n+\beta)(\delta)_{pn,k}}\int\limits_{0}^{1}t^{\frac{a+n{\alpha}}{k}-1}(1-t)^{\frac{b}{k}-1}\exp\left(\frac{-A}{t^{\frac{m}{k}}(1-t)^{\frac{m}{k}}}\right)dt $$
By using (1.14) as in the above equation, we attain the required result.

\noindent{\bf Corollary 2.1.} For $A=0,a={\beta}$ in Theorem 2.1, we deduce the following result:
$$\frac{1}{{\Gamma}_{k}{(b)}}\int\limits_{0}^{1}t^{\frac{\beta}{k}-1}(1-t)^{\frac{b}{k}-1}E_{k,\alpha,\beta,p}^{\gamma,\delta,q}(zt^\frac{\alpha}{k})dt$$
$$=\sum\limits_{n=0}^{\infty}\frac{(\gamma)_{qn,k}{z^{n}}}{(\delta)_{pn,k}{\Gamma}_{k}(\beta+b+{\alpha}n)}.\eqno(2.2)$$

\noindent{\bf Theorem 2.2.} If ${k\in \mathbb R, \alpha,\beta,\gamma,\delta,\rho,\mu} \in \mathbb C, \Re(\alpha)>0, \Re(\beta)>0, \Re(\gamma)> 0, \Re(\delta)> 0, \Re(\rho), \Re(\mu), \Re(A)>0 $, $p,q >0$ and $q < R{\alpha} + p $; $\mid arg(\frac{b^{'}c+d}{a^{'}c+d})\mid<{\pi}$, then

$$\int\limits_{t}^{x}(x-t)^{\frac{\rho}{k}-1}(s-t)^{\frac{\mu}{k}-1}\exp\left(\frac{-A}{(x-s)^{\frac{m}{k}}(s-t)^{\frac{m}{k}}}\right)E_{k,\alpha,\beta, p}^{\gamma,\delta,q}(z(s-t)^{\frac{\alpha}{k}})dt$$

$$=\sum\limits_{r=0}^{\infty}\sum\limits_{n=0}^{\infty}\dfrac{(x-t)^{\frac{\rho+\mu+n{\alpha}-2mr}{k}-2}(-A)^{r}(\gamma)_{qn,k}{z^{n}}}{{\Gamma}_{k}({\alpha}n+\beta)(\delta)_{pn,k}r!}B_{k}(\mu+n{\alpha}-mr,\rho-mr).\eqno(2.3)$$

{\bf Proof.} On taking L.H.S. of (2.2) and then by changing the variable s to $u=\frac{s-t}{x-t}$, we get:

$$\int\limits_{t}^{x}(x-t)^{\frac{\rho}{k}-1}(s-t)^{\frac{\mu}{k}-1}\exp\left(\frac{-A}{(x-s)^{\frac{m}{k}}(s-t)^{\frac{m}{k}}}\right)E_{k,\alpha,\beta, p}^{\gamma,\delta,q}(z(s-t)^{\frac{\alpha}{k}})dt,$$

$$=\int\limits_{0}^{1}(1-u)^{\frac{\rho}{k}-1}(x-t)^{\frac{\rho}{k}-1}u^{\frac{\mu}{k}-1}(x-t)^{\frac{\mu}{k}-1}\exp\left(\frac{-A}{((x-t)(1-u))^{\frac{m}{k}}(u(x-t))^{\frac{m}{k}}}\right)E_{k,\alpha,\beta, p}^{\gamma,\delta,q}(zu^{\frac{\alpha}{k}}(x-t)^{\frac{\alpha}{k}})dt.$$
Expanding the exponential function and k-Mittag-Leffler function in their respective series, we attain:

$$\int\limits_{0}^{1}(1-u)^{\frac{\rho}{k}-1}(x-t)^{\frac{\rho}{k}-1}u^{\frac{\mu}{k}-1}(x-t)^{\frac{\mu}{k}-1}\sum\limits_{r=0}^{\infty}\frac{(-A)^{r}}{(x-t)^{\frac{2mr}{k}}(1-u)^{\frac{mr}{k}}u^{\frac{mr}{k}}r!}\sum\limits_{n=0}^{\infty}\frac{(\gamma)_{qn,k}{z^{n}u^{\frac{n{\alpha}}{k}}(x-t)^{\frac{\alpha{n}}{k}}}}{{\Gamma}_{k}({\alpha}n+\beta)(\delta)_{pn,k}}.$$
By changing the order of summation and integration (which is guaranteed under the given conditions), we get:
$$=\sum\limits_{r=0}^{\infty}\sum\limits_{n=0}^{\infty}\dfrac{(x-t)^{\frac{\rho+\mu+n{\alpha}-2mr}{k}-2}(-A)^{r}(\gamma)_{qn,k}{z^{n}}}{{\Gamma}_{k}({\alpha}n+\beta)(\delta)_{pn,k}r!}\int\limits_{0}^{1}(1-u)^{\frac{\rho-mr}{k}-1}u^{\frac{\mu+n{\alpha-mr}}{k}-1}dt,$$
which further on using the integral(1.11) yields the required result.

\noindent{\bf Corollary 2.2.} For $A=0,\mu=\beta$ in Theorem 2.2, we get:
$$\frac{1}{\Gamma_{k} {\rho}}\int\limits_{t}^{x}(x-t)^{\frac{\rho}{k}-1}(s-t)^{\frac{\beta}{k}-1}E_{k,\alpha,\beta, p}^{\gamma,\delta,q}(z(s-t)^{\frac{\alpha}{k}})dt
=\sum\limits_{n=0}^{\infty}\dfrac{(x-t)^{\frac{\rho+\beta+n{\alpha}}{k}-2}(\gamma)_{qn,k}{z^{n}}}{{\Gamma}_{k}({\alpha}n+\beta)(\delta)_{pn,k}}.\eqno(2.4)$$

\noindent{\bf Theorem 2.3.} If ${k\in \mathbb R, \alpha,\beta,\gamma,\delta,\rho,\mu, \lambda,\sigma} \in \mathbb C, \Re(\alpha)>0, \Re(\beta)>0, \Re(\gamma)>0,\Re(\sigma) > 0, \Re(\delta)> 0, \Re(\rho), \Re(\mu), \Re(\lambda)>0, \Re(A)>0 $; $p,q >0$ and $q < \Re{\alpha} + p $, then
$$\int\limits_{0}^{1}t^{\frac{\lambda}{k}-1}(1-t)^{\frac{\mu-\lambda}{k}-1}\left(1-ut^{\frac{\rho}{k}}(1-t)^{\frac{\sigma}{k}}\right)^{-a}\exp\left(\frac{-A}{t^{\frac{m}{k}}(1-t)^{\frac{m}{k}}}\right)E_{k,\alpha,\beta, p}^{\gamma,\delta,q}(zt^{\frac{\alpha}{k}})dt$$

$$=\sum\limits_{r=0}^{\infty}\sum\limits_{n=0}^{\infty}\dfrac{(a)_{r}u^{r}(\gamma)_{qn,k}{z^{n}}}{{\Gamma}_{k}({\alpha}n+\beta)(\delta)_{pn,k}r!}B_{k}(\lambda+n{\alpha}+{\rho}r,\mu-\lambda+r{\sigma};m;A).\eqno(2.5)$$

\noindent{\bf Proof.} On taking L.H.S. of Theorem 2.3, using the definition of generalized k-Mittag-Leffler function (1.10), and then by changing the order of summation and integration, we get: $$\int\limits_{0}^{1}t^{\frac{\lambda}{k}-1}(1-t)^{\frac{\mu-\lambda}{k}-1}\left(1-ut^{\frac{\rho}{k}}(1-t)^{\frac{\sigma}{k}}\right)^{-a}\exp\left(\frac{-A}{t^{\frac{m}{k}}(1-t)^{\frac{m}{k}}}\right)E_{k,\alpha,\beta, p}^{\gamma,\delta,q}(zt^{\frac{\alpha}{k}})dt $$

$$=\int\limits_{0}^{1}t^{\frac{\lambda}{k}-1}(1-t)^{\frac{\mu-\lambda}{k}-1}\sum\limits_{r=0}^{\infty}\frac{(a)_{r}u^{r}t^{\frac{{\rho}r}{k}}(1-t)^\frac{r{\sigma}}{k}}{r!}\exp\left(\frac{-A}{t^{\frac{m}{k}}(1-t)^{\frac{m}{k}}}\right)\sum\limits_{n=0}^{\infty}\frac{(\gamma)_{qn,k}{z^{n}t^{\frac{n \alpha}{k}}}}{{\Gamma}_{k}({\alpha}n+\beta)(\delta)_{pn,k}}dt,$$

$$=\sum\limits_{r=0}^{\infty}\sum\limits_{n=0}^{\infty}\dfrac{(a)_{r}u^{r}(\gamma)_{qn,k}{z^{n}}}{{\Gamma}_{k}({\alpha}n+\beta)(\delta)_{pn,k}r!}\int\limits_{0}^{1}t^{\frac{\lambda+n{\alpha}+{\rho}r}{k}-1}(1-t)^{\frac{\mu-\lambda+{\sigma}r}{k}-1}\exp\left(\frac{-A}{t^{\frac{m}{k}}(1-t)^{\frac{m}{k}}}\right)dt.$$
By using (1.14) as in the above equation, we derive the required result.

{\bf Corollary 2.3.} For $a=0$ in Theorem 2.3 reduces to the following result: $$\int\limits_{0}^{1}t^{\frac{\lambda}{k}-1}(1-t)^{\frac{\mu-\lambda}{k}-1}\exp\left(\frac{-A}{t^{\frac{m}{k}}(1-t)^{\frac{m}{k}}}\right)E_{k,\alpha,\beta, p}^{\gamma,\delta,q}(zt^{\frac{\alpha}{k}})dt$$
$$=\sum\limits_{n=0}^{\infty}\dfrac{(\gamma)_{qn,k}{z^{n}}}{{\Gamma}_{k}({\alpha}n+\beta)(\delta)_{pn,k}}B_{k}(\lambda+n{\alpha}, \mu-\lambda;A).\eqno(2.6)$$

{\bf Corollary 2.4.} Setting $a=0$, $A=0$ in Theorem 2.3, we deduces the following result:
$$\int\limits_{0}^{1}t^{\frac{\lambda}{k}-1}(1-t)^{\frac{\mu-\lambda}{k}-1}E_{k,\alpha,\beta, p}^{\gamma,\delta,q}(zt^{\frac{\alpha}{k}})dt
=\sum\limits_{n=0}^{\infty}\frac{(\gamma)_{qn,k}{z^{n}}}{{\Gamma}_{k}({\alpha}n+\beta)(\delta)_{pn,k}}B_{k}(\lambda+n{\alpha}, \mu-\lambda).\eqno(2.7)$$

\noindent{\bf Remark.} If we consider $p=q=k=m=1$ in Theorem (2.1), (2.2) and (2.3), we get a new class of Beta type integral operators involving the generalized Mittag-Leffler function defined by Salim [22] and the case $\delta=p=k=1$ of (2.1), (2.3) and (2.5) is seen to yield the known results of Ahmed and Khan [1].

\noindent{\bf 3. Special Cases }
\noindent\par

In this section, we establish the following potentially useful integral operators involving generalized k-Beta type functions as special cases of our main results:\\

\noindent{\bf 1.} On setting $\gamma=q=1$ in Theorem 2.1, we get:
$$\int\limits_{0}^{1}t^{\frac{a}{k}-1}(1-t)^{\frac{b}{k}-1}\exp\left(\frac{-A}{t^{\frac{m}{k}}(1-t)^{\frac{m}{k}}}\right)E_{k,\alpha,\beta, p}^{\delta}(zt^{\frac{\alpha}{k}})dt$$
$$=\sum\limits_{n=0}^{\infty}\dfrac{{z^{n}}}{\Gamma_{k}({\alpha}n+\beta)(\delta)_{pn,k}}B_{k}(n{\alpha}+a,b;A;m).\eqno(3.1)$$

\noindent{\bf 2.} On setting $\alpha=\beta=q=\gamma=1$ in Theorem 2.1, we find:
$$\int\limits_{0}^{1}t^{\frac{a}{k}-1}(1-t)^{\frac{b}{k}-1}\exp\left(\frac{-A}{t^{\frac{m}{k}}(1-t)^{\frac{m}{k}}}\right)E_{k, p}^{\delta}(zt^{\frac{1}{k}})dt$$
$$=\sum\limits_{n=0}^{\infty}\dfrac{{z^{n}}}{(\delta)_{pn,k}}B_{k}(n+a,b;A;m).\eqno(3.2)$$

\noindent{\bf 3.} On setting $\delta=p=1$ in Theorem 2.1, we attain:
$$\int\limits_{0}^{1}t^{\frac{a}{k}-1}(1-t)^{\frac{b}{k}-1}\exp\left(\frac{-A}{t^{\frac{m}{k}}(1-t)^{\frac{m}{k}}}\right)E_{k,\alpha, \beta}^{\gamma,q}(zt^{\frac{\alpha}{k}})dt$$
$$=\sum\limits_{n=0}^{\infty}\dfrac{{(\gamma)_{qn}z^{n}}}{\Gamma_{k}({\alpha}n+\beta)}B_{k}(n{\alpha}+a,b;A;m).\eqno(3.3)$$

\noindent{\bf 4.} On setting $\gamma=q=1$ in Theorem 2.2, we achieve:
$$\int\limits_{t}^{x}(x-t)^{\frac{\rho}{k}-1}(s-t)^{\frac{\mu}{k}-1}\exp\left(\frac{-A}{(x-s)^{\frac{m}{k}}(s-t)^{\frac{m}{k}}}\right)E_{k,\alpha,\beta, p}^{\delta}(zt^{\frac{\alpha}{k}})dt$$
$$=\sum\limits_{r=0}^{\infty}\sum\limits_{n=0}^{\infty}\dfrac{(x-t)^{\frac{\rho+\mu+n{\alpha}-2mr}{k}-2}(-A)^{r}{z^{n}}}{{\Gamma}_{k}({\alpha}n+\beta)(\delta)_{pn,k}r!}B_{k}(\mu+n{\alpha}-mr,\rho-mr).\eqno(3.4)$$

\noindent{\bf 5.} On setting ${\alpha}= {\beta}= {\gamma} = q=1$ in Theorem 2.2, we acquire:
$$\int\limits_{t}^{x}(x-t)^{\frac{\rho}{k}-1}(s-t)^{\frac{\mu}{k}-1}\exp\left(\frac{-A}{(x-s)^{\frac{m}{k}}(s-t)^{\frac{m}{k}}}\right)E_{k, p}^{\delta}(zt^{\frac{1}{k}})dt$$
$$=\sum\limits_{r=0}^{\infty}\sum\limits_{n=0}^{\infty}\dfrac{(x-t)^{\frac{\rho+\mu+n-2mr}{k}-2}(-A)^{r}{z^{n}}}{(\delta)_{pn,k}r!}B_{k}(\mu+n-mr,\rho-mr).\eqno(3.5)$$

\noindent{\bf 6.} On setting ${\delta}= p=1$ in Theorem 2.2, we found:
$$\int\limits_{t}^{x}(x-t)^{\frac{\rho}{k}-1}(s-t)^{\frac{\mu}{k}-1}\exp\left(\frac{-A}{(x-s)^{\frac{m}{k}}(s-t)^{\frac{m}{k}}}\right)E_{k,\alpha,\beta}^{\gamma,q}(zt^{\frac{\alpha}{k}})dt$$
$$=\sum\limits_{r=0}^{\infty}\sum\limits_{n=0}^{\infty}\dfrac{(x-t)^{\frac{\rho+\mu+n{\alpha}-2mr}{k}-2}(-A)^{r}{z^{n}}}{(\delta)_{pn,k}r!}B_{k}(\mu+n{\alpha}-mr,\rho-mr).\eqno(3.6)$$

\noindent{\bf 7.} On setting ${\gamma= q = 1}$ in Theorem 2.3, we find:
$$\int\limits_{0}^{1}t^{\frac{\lambda}{k}-1}(1-t)^{\frac{\mu-\lambda}{k}-1}\left(1-ut^{\frac{\rho}{k}}(1-t)^{\frac{\sigma}{k}}\right)^{-a}\exp\left(\frac{-A}{t^{\frac{m}{k}}(1-t)^{\frac{m}{k}}}\right)E_{k,\alpha,\beta, p}^{\delta}(zt^{\frac{\alpha}{k}})dt $$
$$=\sum\limits_{r=0}^{\infty}\sum\limits_{n=0}^{\infty}\dfrac{(a)_{r}u^{r}{z^{n}}}{{\Gamma}_{k}({\alpha}n+\beta)(\delta)_{pn,k}r!}B_{k}(\lambda+n{\alpha}+{\rho}r,\mu-\lambda+r{\sigma};m;A).\eqno(3.7)$$

\noindent{\bf 8.} On setting ${\alpha}= {\beta}={\gamma}= q=1$ in Theorem 2.3, we get:
$$\int\limits_{0}^{1}t^{\frac{\lambda}{k}-1}(1-t)^{\frac{\mu-\lambda}{k}-1}\left(1-ut^{\frac{\rho}{k}}(1-t)^{\frac{\sigma}{k}}\right)^{-a}\exp\left(\frac{-A}{t^{\frac{m}{k}}(1-t)^{\frac{m}{k}}}\right)E_{k,p}^{\delta}(zt^{\frac{1}{k}})dt $$
$$=\sum\limits_{r=0}^{\infty}\sum\limits_{n=0}^{\infty}\dfrac{(a)_{r}u^{r}{z^{n}}}{(\delta)_{pn}{r!}}B_{k}(\lambda+n+{\rho}r,\mu-\lambda+r{\sigma};m;A).\eqno(3.8)$$

\vspace{.35cm}
\parindent=0mm
{\bf References}\\
\begin{enumerate}
\item[{[1]}] Ahmed, S, Khan, M. A, Euler type integral involving generalization Mittag-Leffler function, Commun. Korean Math. Soc., {\bf{29}}(3)(2014), 479-487.
\item[{[2]}] Chaudhry, M. A, Qadir, A, Rafiq, M and Zubair, S. M, Extension of Euler's beta function, J. Comput. Appl. Math., {\bf{78}}(1)(1997),19-32.
\item[{[3]}] Chaudhry, M. A, Qadir, A and Srivastava, H. M, and Paris R. B, Extended hypergeometric and confluent hypergeometric functions, Appl. Math. Comput., {\bf{159}}(2)(2004), 589-602.
\item [{[4]}] Chaudhry, M. A and Zubair, S. M, Generalized incomplete gamma functions with applications, J. Comput. Appl. Math., {\bf{75}}(1)(1994), 99-124.
\item[{[5]}] Dorrego, G. A and Cerutti, R. A, The k-Mittag-Leffler function, Int. J. Contemp. Math. Sciences, {\bf {15}}(7)(2012), 705-716.
\item[{[6]}] Diaz, R, and Pariguan, E, On hypergeometric functions and k-Pochhammer symbol, Divulgaciones Matematicas,  {\bf {15}}(2)(2007), 179-192.
\item[{[7]}] Gehlot, K. S and Prajapati, J. C, Fractional calculus of generalized k-Wright function, J. Fract. Calc. Appl., {\bf {4}}(2)(2013), 283-289.
\item[{[8]}]Gupta, A and Parihar, C. L, K-new generalized Mittag-Leffler function, J. Frac. Cal. Appl., {\bf {5}}{(1)(2014), 165-176.}
\item[{[9]}] Kilbas, A. A, Saigo, M, Saxena, R. K, Generalized Mittag-leffler function and generalized fractional calculus operators, Int. Tran. Spec. Funct., {\bf {15}}(2004), 31-49.
\item[{[10]}] Kilbas, A, Srivastava, H.M and Trujillo, J, Theory and applications of fractional differential equations, Elsevier, 2006.
\item[{[11]}] Kokologiannaki, Ch, Properties and inequalities of generalized k-gamma,beta and zeta functions, Int. J. Contemp. Math. Science, {\bf {5}}(2010, 653-660.
\item[{[12]}] Krasniqi, V, A limit for the k-gamma and k-beta function, Int. Math. Forum, {\bf{33}}(5)(2010), 1613-1617.
\item[{[13]}] Lee, D. M, Rathie, A. K, Parmar, R. K, and  Kim, Y. S. Generalization of extended beta function, Hypergeometric and confluent Hypergeometric function, Honam Mathematical Journal, {\bf {33}}(3)(2011), 187-206.
\item [{[14]}] Luke, Y. L, The special functions and their approximations, Vol.1 New York, Academic Press, 1969.
\item [{[15]}] Merovci, F, Power product inequalities for the $\Gamma_{k}$ function, Int. J. of Math. Analysis, {\bf {4}}(21)(2010), 1007-1012.
\item [{[16]}] Mittag-Leffler, G.M, and Sur la nouvelle function $E_{\alpha}(x)$, C. R .Acad. Sci Paris, {\bf{137}}(1903), 554-558.
\item[{[17]}] Mansour, M, Determinig the k-generalized gamma function $\Gamma(x)$ by functional equations. Int. J. Contemp. Math. Science, {\bf{21}}(4)(2009), 1037-1042.
\item [{[18]}] Prabhakar, T.R, A singular integral equation with a generalized Mittag-Leffler function in the kernal, Yokohama Math. J., {\bf {19}}(1971), 7-15.
\item[{[19]}] Prudnikov, A.P, Brychkov, Yu.A and Matichev, O. I, Integrals and Series, Vol. I, Gordan and Breach Science Publishers, New York, 1990.
\item[{[20]}] Rainville, E. D, Special functions, the Macmillan Company, New York, 1913.
\item[{[21]}] Salim, T. O, and Faraj, W, A generalization of Mittag-Leffler function and integral operator associated with fractional calculus, Appl. Math. Comput., {\bf {3}}(5)(2012),1-13.
\item[{[22]}] Salim, T. O, Some properties relating to the generalized Mittag-Leffler function, Adv. Appl. Math. Anal., {\bf {4}}(2009), 21-80.
\item[{[23]}] Soubhia, A, Camargo, R and De Oliveira, E. J. Vaz. Theorem for series in the three parameter Mittag-Leffler function, Frac. Calc. Appl. Anal., {\bf {13}}(1)(2010).
\item[{[24]}] Shukla, A.K, Prajapati, J.C, On a generalization of Mittag-Leffler function and its properties, J. Math. Ann. Appl., {\bf{336}}(2007),
797-811.
\item[{[25]}] Srivastava, H.M and Tomovaski, Z, Fractional calculus with an integral operator containing a generalized Mittag-leffler function in the kernel, Appl. Math. Comput., {\bf {211}}(2009), 198-210.
\item[{[26]}] Wiman, A, Uber den fundamental satz in der theory der funktionen, Acta Math., {\bf{29}}(1905), 191-201.
\item[{[27]}] Wiman, A, Uber die Nullstellun der  funktionen, Acta Math., {\bf{29}}(1905), 217-234.
\end{enumerate}

\end{document}